\documentclass[reqno,centertags, 12pt]{amsart}
\usepackage{amsmath,amsthm,amscd,amssymb}
\usepackage{latexsym}
\newcommand{\bbA}{{\mathbb{A}}}
\newcommand{\bbC}{{\mathbb{C}}}
\newcommand{\bbD}{{\mathbb{D}}}

\newcommand{\bbR}{{\mathbb{R}}}

\newcommand{\lb}{\label}
\newcommand{\f}{\frac}
\newcommand{\ul}{\underline}
\newcommand{\ol}{\overline}
\newcommand{\ti}{\tilde  }

\newcommand{\spec}{\text{\rm{spec}}}

\newcommand{\s}{\text{\rm{s}}}

\newcommand{\supp}{\text{\rm{supp}}}

\newcommand{\bi}{\bibitem}

\newcommand{\beq}{\begin{equation}}
\newcommand{\eeq}{\end{equation}}
\newcommand{\ba}{\begin{align}}
\newcommand{\ea}{\end{align}}
\newcommand{\veps}{\varepsilon}

\newcounter{smalllist}
\newenvironment{SL}{\begin{list}{{\rm\roman{smalllist})}}{%
\setlength{\topsep}{0mm}\setlength{\parsep}{0mm}\setlength{\itemsep}{0mm}%
\setlength{\labelwidth}{2em}\setlength{\leftmargin}{2em}\usecounter{smalllist}%
}}{\end{list}}

\DeclareMathOperator{\Sz}{Sz}
\DeclareMathOperator{\Real}{Re}
\DeclareMathOperator{\Ima}{Im}

\allowdisplaybreaks
\numberwithin{equation}{section}

\newtheorem{theorem}{Theorem}[section]
\newtheorem{proposition}[theorem]{Proposition}

\theoremstyle{definition}

\theoremstyle{remark}

\newcommand{\abs}[1]{\lvert#1\rvert}

\begin{document}

\title[Decay and Analyticity]{Jost Functions and Jost Solutions for Jacobi Matrices,
II.~Decay and Analyticity}
\author[D. Damanik and B. Simon]{David Damanik$^{1,2}$ and Barry Simon$^{1,3}$}

\thanks{$^1$ Mathematics 253-37, California Institute of Technology, Pasadena, CA 91125.
E-mail: damanik@caltech.edu; bsimon@caltech.edu}
\thanks{$^2$ Supported in part by NSF grant DMS-0227089}
\thanks{$^3$ Supported in part by NSF grant DMS-0140592}

\date{January 5, 2005}

\begin{abstract} We present necessary and sufficient conditions on the
Jost function for the corresponding Jacobi parameters $a_n -1$ and $b_n$
to have a given degree of exponential decay.
\end{abstract}

\maketitle

\section{Introduction} \lb{s1}

Among the most interesting results in spectral theory are those that give equivalent
sets of conditions --- one set involving recursion coefficients and the other
involving spectral data. Examples are Verblunsky's version \cite{V36} of the Szeg\H{o}
theorem (see \cite{OPUC1}), the strong Szeg\H{o} theorem written as a sum rule
(see \cite{OPUC1}), the Killip-Simon theorem \cite{KS} characterizing $L^2$
perturbations of the free Jacobi matrix, and Baxter's theorem \cite{Bax,OPUC1}.

Our goal in this paper is to present such an equivalence for Jacobi matrices
concerning exponential decay. That is, we consider orthogonal polynomials on
the real line (OPRL) whose recursion relation is
\begin{equation} \lb{1.1}
xp_n(x) = a_{n+1} p_{n+1}(x) + b_{n+1} p_n(x) + a_n p_{n-1}(x)
\end{equation}
for Jacobi parameters $\{a_n\}_{n=1}^\infty$, $\{b_n\}_{n=1}^\infty$. Here $p_n(x)$
are the orthonormal polynomials and $p_{-1}(x)\equiv 0$ (i.e., $a_0$ is not needed in
\eqref{1.1} for $n=0$).

\eqref{1.1} is often summarized by the Jacobi matrix
\begin{equation} \lb{1.1a}
J=\begin{pmatrix} b_1 & a_1 & 0 & \dots \\
a_1 & b_2 & a_2 & \dots \\
0 & a_2 & b_3 & \dots \\
\dots & \dots & \dots & \dots \\
\dots & \dots & \dots & \dots
\end{pmatrix}
\end{equation}
By $J_0$ we mean the $J$ with $a_n \equiv 1$, $b_n\equiv 0$.

The model of what we will find here is the following result of Nevai-Totik \cite{NT89}
in the theory of orthogonal polynomials on the unit circle (OPUC):

\begin{theorem}[Nevai-Totik \cite{NT89}; see Section~7.1 of \cite{OPUC1}]\lb{T1.1}
Let $d\mu$ be a probability measure on $\partial\bbD$ obeying
\begin{equation} \lb{1.2}
d\mu = w(\theta)\, \f{d\theta}{2\pi} + d\mu_\s
\end{equation}
Fix $R>1$. Then the following are equivalent:
\begin{SL}
\item[{\rm{(1)}}] The Szeg\H{o} condition holds, $d\mu_\s=0$, and the Szeg\H{o}
function, $D(z)$, has $D(z)^{-1}$ analytic in $\{z\mid \abs{z}<R\}$.
\item[{\rm{(2)}}]
\begin{equation} \lb{1.3}
\limsup_{n\to\infty}\, \abs{\alpha_n}^{1/n} \leq R^{-1}
\end{equation}
\end{SL}
\end{theorem}

The Szeg\H{o} condition is
\begin{equation} \lb{1.4}
\int \log(w(\theta))\, \f{d\theta}{2\pi} > -\infty
\end{equation}
in which case $D$ is defined initially on $\bbD$ by
\begin{equation} \lb{1.5}
D(z) =\exp\biggl( \int \f{e^{i\theta}+z}{e^{i\theta}-z}\, \log(w(\theta))\,
\f{d\theta}{4\pi}\biggr)
\end{equation}
In \eqref{1.3}, $\alpha_n$ are the Verblunsky coefficients, that is, the recursion
coefficients for the monic OPUC, $\Phi_n$,
\begin{equation} \lb{1.6}
\Phi_{n+1}(z) = z\Phi_n(z) -\bar\alpha_n \Phi_n^*(z)
\end{equation}
with
\begin{equation} \lb{1.7}
\Phi_n^*(z) = z^n\, \ol{\Phi_n (1/\bar z)}
\end{equation}
See \cite{OPUC1,OPUC2,Szb,GBk,GBk1} for background on OPUC.

Also relevant to our motivation is the following simple result:

\begin{theorem}\lb{T1.2} Let $d\mu$ be a probability measure on $\partial\bbD$
obeying \eqref{1.2}. Then the following are equivalent:
\begin{SL}
\item[{\rm{(1)}}] The Szeg\H{o} condition holds, $d\mu_\s =0$, and the Szeg\H{o}
function, $D(z)$, has $D(z)^{-1}$ a polynomial of exact degree $n$.
\item[{\rm{(2)}}] $\alpha_j =0$ for $j\geq n$ and $\alpha_{n-1}\neq 0$.
\end{SL}
\end{theorem}

\begin{proof} {\ul{(2) $\Rightarrow$ (1).}} \ In this case (see
\cite[Theorem~1.7.8]{OPUC1}), $d\mu = \f{d\theta}{2\pi}
\abs{\varphi_n^* (e^{i\theta})}^{-2}$, so $D^{-1} =\varphi_n^*$ is
a polynomial.

\smallskip
\noindent {\ul{(1) $\Rightarrow$ (2).}} \ $D(z)^{-1}$ is nonvanishing on
$\overline{\bbD}$, and so the measure has the form $\f{d\theta}{2\pi}\abs{D(z)}^2$ and so
has $\alpha_j =0$ for $j \geq n$ (\cite[Theorem~1.7.8]{OPUC1}). Thus $D(z)^{-1} =
\varphi_n^*(z)$, and since $\varphi_n$ has degree exactly $n$, $\Phi_n^* = \Phi_{n-1}^*
-\alpha_{n-1} z\Phi_{n-1}$ implies $\alpha_{n-1} \neq 0$.
\end{proof}

In our work, the spectral measure has the form
\begin{equation} \lb{1.8}
d\gamma(x) =f(x)\, dx + d\gamma_\s
\end{equation}
where $\supp f \subset [-2,2]$. We say $d\gamma_\s$ is regular if $\gamma_\s
([-2,2])=0$ and $d\gamma_\s$ has finite support (i.e., no embedded singular
spectrum and only finitely many bound states). The $m$-function associated to
$d\gamma$ is defined on $\bbC\backslash\supp(d\gamma)$ by
\begin{equation} \lb{1.9}
m(E) = \int \f{d\gamma(x)}{x-E}
\end{equation}
and $M$ is defined on $\bbD =\{z\mid\abs{z} <1\}$ by
\begin{equation} \lb{1.10}
M(z) = -m(z+z^{-1})
\end{equation}
Since $z\mapsto z+z^{-1}$ maps $\bbD$ to $\bbC\cup\{\infty\}\backslash [-2,2]$,
$M$ is analytic on $\bbD\backslash \{z\in\bbR\cap\bbD \mid z+z^{-1}$ is a
point mass of $d\gamma\}$ with simple poles at the missing points.

The Jost function, $u(z)$, is defined and analytic on $\bbD$ in many cases
and determined first by
\begin{equation} \lb{1.11}
\abs{u(e^{i\theta})}^2 \Ima M(e^{i\theta}) =\sin\theta
\end{equation}
where the functions at $e^{i\theta}$ are a.e.\ limits as $r\uparrow 1$ of
the functions at $re^{i\theta}$. The second condition on $u$ is that, for
$z\in\bbD$,
\begin{equation} \lb{1.12}
u(z) =0 \Leftrightarrow z+z^{-1}\text{ is a point mass of } d\gamma
\end{equation}

If one has the sufficient regularity of $\Ima M$ on $\partial\bbD$ and
$\gamma_\s$ is regular, \eqref{1.11}/\eqref{1.12} determine
$u$ via
\begin{equation} \lb{1.13x}
u(z) = \prod_{u(z_j)=0} \biggl( \f{z-z_j}{1-\bar z_j z}\biggr) \exp \biggl( \int
\f{e^{i\theta}+z}{e^{i\theta}-z}\, \log \biggl( \f{\sin\theta}{\Ima M(\theta)}\biggr)
\f{d\theta}{4\pi}\biggr)
\end{equation}
In addition, if the Jacobi parameters obey
\[
\sum_{n=1}^\infty \, \abs{a_n -1} + \abs{b_n} <\infty
\]
then the Jost function can be directly constructed using variation of parameters (see
Teschl \cite{Teschl}), perturbation determinants (see Killip-Simon \cite{KS}), or an
approach of Geronimo-Case \cite{GC80}. Since this latter approach is not well-known and
those authors do not provide the detailed estimates we will need, we have described this
approach in Appendix~A.

When there are zeros of $u$ in $\bbD$, then $u$ does not uniquely determine
$d\gamma$. $f$ is determined by \eqref{1.11} and
\begin{equation} \lb{1.12x}
f(2\cos\theta) =\Ima M(e^{i\theta})
\end{equation}
and the positions of the point masses are the zeros, but the weights, $w_j$, of the zeros
(i.e., the values of $\gamma (\{E_j\})=w_j$) are needed. The possible values of $w_j$ are
constrained by
\begin{equation} \lb{1.13}
\sum_j w_j + 2\int_0^\pi \f{\sin^2\theta}{\abs{u(e^{i\theta})}^2}\, d\theta =1
\end{equation}
by \eqref{1.12x}, \eqref{1.11}, and
\begin{equation} \lb{1.14}
\int_{-2}^2 f(E)\, dE = 2\int_0^\pi f(2\cos\theta) \sin\theta\, d\theta
\end{equation}

Thus, modulo some regularity issues, the knowledge of a $d\gamma$ with regular
$d\gamma_\s$ is equivalent to the knowledge of $u$ and the finite number of
weights $w_j$ constrained by \eqref{1.13}. Our main goal in this paper is to
describe what Jost functions and weights are associated to $a_n$'s and $b_n$'s
with a given rate of exponential decay or with finite support. We will view the
Jost function/weights as spectral data. This is justified by the following:

\begin{theorem}\lb{T1.3} Let $u$ be a function analytic in a neighborhood of
$\bar\bbD$ whose only zeros in this neighborhood lie in $\overline{\bbD} \cap \bbR$ with
those zeros all simple. For each zero in $\bbD\cap\bbR$, let a weight $w_j >0$ be given
so that \eqref{1.13} holds. Then there is a unique measure $d\gamma$ for which $u$ is the
Jost function and $w_j$ the weights.
\end{theorem}

Since this is peripheral to the main thrust of this paper, we do not give a
detailed proof, but note several remarks:
\begin{SL}
\item[1.] Related issues are discussed in Paper~I of this series \cite{Jost1}. \item[2.]
One first shows that the $M$ defined by $d\gamma$ has a meromorphic continuation to a
neighborhood of $\overline{\bbD}$; this is done in Theorem~13.7.1 of \cite{OPUC2}.
\item[3.] The methods we use in Sections~\ref{s2} and \ref{s3} then show that the $a_n
-1$ and $b_n$ decay exponentially. \item[4.] Thus, by the results of Appendix~A, a Jost
function, $\ti u$, exists. $u/\ti u$ has removable singularities, is nonzero on $\bbD$,
is analytic in a neighborhood of $\overline{\bbD}$, and on $\partial\bbD$, $\abs{u/\ti
u}=1$. Thus, $u=\bar u$.
\end{SL}

\smallskip
The perturbation determinant can be defined by
\begin{equation} \lb{1.15}
L(z) = \f{u(z)}{u(0)}
\end{equation}
This is obviously normalized by
\begin{equation} \lb{1.16}
L(0) =1
\end{equation}
which is simpler than \eqref{1.13}. Of course, $u(0)$ can be recovered from
$\{w_j\}_{j=1}^N$ and $L(z)$ by \eqref{1.15} and \eqref{1.13}. We note that
when $J-J_0$ is trace class, we have (see \cite{KS})
\begin{equation} \lb{1.17}
L(z) =\det (1+(J-J_0)[J_0 - (z+z^{-1})]^{-1})
\end{equation}

Our goal in this paper is to prove four theorems: two in the simple case
where there is no point spectrum and two in the general case. In each pair,
one describes finite support perturbations and one, exponential decay.
We begin with the case of no bound states:

\begin{theorem}\lb{T1.3A} If $a_n =1$ and $b_n =0$ for large $n$, then
$L(z)$ is a polynomial. Conversely, any polynomial $L(z)$ which obeys
\begin{SL}
\item[{\rm{(i)}}] $L(z)$ is nonvanishing on $\overline{\bbD}\backslash\{\pm 1\}$
\item[{\rm{(ii)}}] If $+1$ and/or $-1$ are zeros, they are simple \item[{\rm{(iii)}}]
$L(0)=1$
\end{SL}
is the perturbation determinant of a unique Jacobi matrix and it obeys
$a_n =1$ and $b_n =0$ for all large $n$.
\end{theorem}

{\it Remark.} By Theorem~\ref{TA.1.1}, there is a precise relation between the degree of
$L$ and the range of $(a_n - 1, b_n)$.

\begin{theorem} \lb{T1.4} Let $R>1$. If
\begin{equation} \lb{1.18}
\lim_{n\to\infty}\, (\abs{a_n -1} + \abs{b_n})^{1/2n} \leq R^{-1}
\end{equation}
then $L(z)$ has an analytic continuation to $\{z\mid\abs{z}<R\}$.
Conversely, if $L(z)$ is analytic in $\{z\mid\abs{z}<R\}$ and obeys
{\rm{(i)--(iii)}} from Theorem~\ref{T1.3A}, then \eqref{1.18} holds.
\end{theorem}

For bound states, things are more complicated. One way of understanding this is the
following. Consider a finite support set of Jacobi parameters with $\ell$ bound states.
One can change the $\ell$ weights without changing $L(z)$. It is known in that case that
changing a single weight introduces an explicit (in terms of solutions of the original
Jacobi recursion) correction which decays exponentially (see \cite[pp.~64-66]{CS} and
\cite{FGJFA,FGGT}). That means one expects only one out of the $\ell$-parameter family of
Jacobi matrices to have Jacobi parameters with finite support. Put differently, there
will be many Jacobi matrices with bound states but only exponentially decaying Jacobi
parameters that have polynomial Jost functions. So the key is identifying the weights
that single out finite support.

Rather than discuss weights, it is more convenient to use residues of poles
of $M$\!. Of course, if $z_j\in\bbD$, $z_j + z_j^{-1} =E_j$ is a point mass
in $\gamma$, then
\begin{equation} \lb{1.19}
w_j = \lim_{E\to E_j}\, (E-E_j) m(E) = (z_j^{-1} -z_j) z_j^{-1} \lim_{z\to z_j}\, (z-z_j)
M(z)
\end{equation}
so the data are equivalent.

\smallskip
\noindent{\bf Definition.} Let $M$ be the $M$-function associated to a $u$
and a set of weights. Suppose $u$ is analytic in $\{z\mid\abs{z}<R\}$ for
some $R>1$ and $u(z_j)=0$ with $\abs{z_j}>R^{-1}$. We say the weight at $z_j$
is canonical if and only if
\begin{equation} \lb{1.20}
\ti w_j \equiv \lim_{z\to z_j}\, (z-z_j) M(z) = - (z_j -z_j^{-1}) [u'(z_j)\, \ol{u(1/\bar
z_j)}\,]^{-1}
\end{equation}

\smallskip
Here are our main theorems on the general case:

\begin{theorem}\lb{T1.5} If $a_n =1$ and $b_n =0$ for large $n$, then $L(z)$
is a polynomial and all the weights are canonical. Conversely, if $L$ is a
polynomial obeying
\[
(\text{\rm{i}}^\prime) \qquad\qquad L(z) \text{ is nonvanishing on }
\overline{\bbD}\backslash\bbR
\]
and {\rm{(ii)--(iii)}} of Theorem~\ref{T1.3A}, then there is at most one set of
Jacobi parameters with $a_n =1$ and $b_n =0$ for $n$ large that has that $L$
as perturbation determinant. Moreover, the weights associated to this set
are the canonical ones. If these canonical weights lead to $w_j >0$, then there
is a set of Jacobi parameters with $a_n =1$ and $b_n =0$ for large $n$.
\end{theorem}

{\it Remark.} It is easy to construct polynomial $L$'s which are not the perturbation
determinant of any finite support Jacobi parameters, although they are perturbation
determinants. For example, if $L(z_0) =L(z_0^{-1}) =0$ for some $z_0 \in (0,1)$,
\eqref{1.20} cannot hold. Thus
\[
L(z) = (1-2z)(1-\tfrac12\,z)
\]
is a perturbation determinant but not for a Jacobi matrix of finite support. There
are also examples where the canonical weights are negative.

\begin{theorem}\lb{T1.6} Let $R>1$. If \eqref{1.18} holds, then $L(z)$ has an
analytic continuation to $\{z\mid\abs{z}<R\}$ and the weights of all $z_j$
with $\abs{z_j}>R^{-1}$ are canonical. Conversely, if $L(z)$ is a function
analytic in $\{z\mid\abs{z}<R\}$ obeying {\rm{(i$^\prime$)--(iii)}} of
Theorem~\ref{T1.5}, then \eqref{1.18} holds if and only if all weights for
$z_j$ with $\abs{z_j}>R^{-1}$ are canonical.
\end{theorem}

These four theorems have a direct part (i.e., going from $\{a_n, b_n\}_{n=1}^\infty$
to $L(z)$) and an inverse part. The direct parts (except for the importance of
canonical weights) are well-known. We provide a proof of all but the canonical
weights in Appendix~A. The canonical weight result is proven in Section~\ref{s3}.

The inverse parts are more subtle --- and the main content of this paper. The
no bound state results appear in Section~\ref{s2} and the bound state results
in Section~\ref{s3}. Our approach is based on the use of coefficient stripping,
that is, relating $u,M$ for $\{a_n,b_n\}_{n=1}^\infty$ to $u,M$ for
$\{\ti a_n, \ti b_n\}_{n=1}^\infty$ where $\ti a_n=a_{n+1}$, $\ti b_n =
b_{n+1}$. Section~\ref{s2} will rely on a remarkably simple contraction argument,
Section~\ref{s3} on the fact that coefficient stripping only preserves
analyticity if weights are canonical.

One can wonder if one can't at least prove the no bound state results by appealing
to the Nevai-Totik theory and the Szeg\H{o} mapping (see Section~13.1 of \cite{OPUC2})
relating OPUC and OPRL. Indeed, we will show in Section~\ref{s2} that our method can
be used to prove the inverse part of their result. There is a difficulty with blind
use of the Szeg\H{o} map, already seen by the fact that $J_0$ does not map into
Verblunsky coefficients with exponential decay (see Example~13.1.3 of \cite{OPUC2}).
This can be understood by noting that the Jost function, $u$, for $d\gamma$ and
the Szeg\H{o} function, $D$, for $\mu= \Sz^{-1}(d\gamma)$ are related by
\begin{equation} \lb{1.21}
D(z)^{-1} = \f{2^{-1/2} u(z)}{1-z^2}
\end{equation}
Thus, $D(z)^{-1}$ is not analytic where $u$ is, unless $u(+1)=u(-1)=0$. In
that case, one can use Nevai-Totik to obtain Theorems~\ref{T1.3A} and \ref{T1.4}.

There are two strategies for dealing with the general case. First (and our
original proof), one can add extra $a$'s and $b$'s at the start to produce
$u(+1) = u(-1)=0$. Second, $\mu\mapsto\Sz(\mu)$ is one of four maps (see
Section~13.2 of \cite{OPUC2}). $\Sz_2$ maps onto all Jacobi matrices
with spectrum on $[-2,2]$ and with $u(1)\neq 0\neq u(-1)$ and has no division
factor. $\Sz_3$ and $\Sz_4$ divide by $1-z$ and $1+z$ and are onto all matrices
with $u(-1)\neq 0$ and $u(1)\neq 0$. In this way, one can always find a $\mu$
with $\gamma =\Sz_j(\mu)$ so the $D$-function for $\mu$ is analytic.

It should also be possible to prove the inverse results we need using the
Marchenko equation. That said, we prefer the approach in Section~\ref{s2}.

Surprisingly, the four main results of this paper appear to be new, although for
Schr\"odinger operators with Yukawa potentials, there are related results in Newton
\cite{New} and Chadan-Sabatier \cite{CS}. Geronimo \cite{Ge94} has a paper closely
related to our theme here, but he makes an a priori hypothesis about $M$ that means his
results are not strictly Jacobi-parameter hypotheses on one side. So he does not have our
results, although it is possible that one can modify his methods to prove them.

An analog of our results on what are Jost functions for Jacobi matrices of finitely
supported Jacobi parameters is the study of the sets of allowed resonance positions for
half-line Schr\"odinger operators with compactly supported potentials. There is a large
literature on this question \cite{Fad,Fro,Koro1,Koro2,Sim,Zwo1,Zwo2}. In particular in
\cite{Koro1,Koro2}, Korotyaev makes some progress in classifying all Jost functions in
this case.

We announced the results in \cite{DKS} and some of them have been presented
in \cite{OPUC2}, but we note an error in \cite{OPUC2}: Theorem~13.7.4 is
wrong because, when stating existence of a finite-range solution, it fails
to require $u(z_j^{-1})\neq 0$ and that the canonical weights be positive.

\section{The Case of No Bound States} \lb{s2}

Our goal in this section is to prove Theorems~\ref{T1.3A} and \ref{T1.4}.
We suppose we have a set of Jacobi parameters $\{a_n,b_n\}_{n=1}^\infty$ with
Jost function, $u(z)\equiv u^{(0)}(z)$, and $M$-function, $M(z)\equiv M^{(0)}(z)$.

Associated to Jacobi parameters $\{a_{k+n}, b_{k+n}\}_{k=1}^\infty$, we have
corresponding Jost function, $u^{(n)}(z)$, and $M$-function, $M^{(n)}(z)$. $u^{(n)}(z)$
is the solution of a difference equation at $0$ where the solution is asymptotic to $z^n$
as $n\to\infty$. It follows that
\begin{equation} \lb{2.1}
u_n(z) = a_n^{-1} z^n u^{(n)}(z)
\end{equation}
obeys (see \eqref{A.1.37})
\begin{equation} \lb{2.2}
a_n u_{n+1} + (b_n -(z+z^{-1}))u_n +a_{n-1} u_{n-1} =0
\end{equation}
Moreover (see \eqref{A.1.36} and \eqref{A.1.38}),
\begin{equation} \lb{2.3}
M^{(n)}(z) = \f{u_{n+1}(z)}{a_n u_n(z)}
\end{equation}

This leads to the following set of update formulae:
\begin{align}
u^{(n+1)}(z) &= a_{n+1} z^{-1} u^{(n)}(z) M^{(n)}(z) \lb{2.4}  \\
M^{(n)}(z)^{-1} &= z+z^{-1} -b_{n+1} -a_{n+1}^2 M^{(n+1)}(z) \lb{2.5}
\end{align}
Since $M(z)=\langle\delta_0, (z+z^{-1} - J)^{-1} \delta_0\rangle$, we see
\begin{equation} \lb{2.6}
\f{M^{(n)}(z)}{z} = 1+ O(z)
\end{equation}
so that \eqref{2.5} implies
\begin{equation} \lb{2.7}
\biggl( \f{M^{(n)}(z)}{z}\biggr)^{-1} = 1-b_{n+1} z - (a_{n+1}^2 -1)
z^2 + O(z^3)
\end{equation}
which means
\begin{equation} \lb{2.8}
\log \biggl( \f{M^{(n)}(z)}{z}\biggr) = b_{n+1} z +
((a_{n+1}^2 -1) + \tfrac12\, b_{n+1}^2)z^2 + O(z^3)
\end{equation}

There is an additional feature we will need. Suppose $u(z)$ is analytic in
$\{z\mid\abs{z}<R\}$ for some $R>1$. Define
\begin{equation} \lb{2.9}
f^\sharp (z) = \ol{f(1/\bar z)}
\end{equation}
for $z\in\bbA_R =\{z\mid R^{-1} <\abs{z} <R\}$. Then we claim
\begin{equation} \lb{2.10}
M(z) -M^\sharp(z) = [u(z) u^\sharp(z)]^{-1} (z-z^{-1})
\end{equation}
To see this, we note (see \eqref{1.11})
\begin{equation} \lb{2.11}
\Ima M(e^{i\theta})= [u(e^{i\theta})\, \ol{u(e^{i\theta})}\,]^{-1} \sin \theta
\end{equation}
This is \eqref{2.10} for $z=e^{i\theta}$, so \eqref{2.10} follows by analyticity.

The strategy of our proof will be to control $u^{(n)}, M^{(n)}$ inductively using
\eqref{2.4}, \eqref{2.5} for $z\in\bbD$ and \eqref{2.10} outside $\bbD$. We will
then feed this control into \eqref{2.8} to control $a_{n+1}-1$ and $b_{n+1}$.

We want to use the update equations to confirm that $u^{(n)}$ is analytic in at
least as big a region as $u$. This will need the assumption that $u$ is
nonvanishing on $\bbD$ (and will be the key issue to be addressed in the next
section).

\begin{theorem}\lb{T2.1} If $u$ is analytic in $\{z\mid\abs{z}<R\}$ and
nonvanishing on $\overline{\bbD}\backslash\{+1, -1\}$ with at most simple zeros at $\pm
1$, then the same is true of each $u^{(n)}$.
\end{theorem}

\begin{proof} By induction, we only need this for $u^{(1)}$. By \eqref{2.10},
$M$ is meromorphic on $\{z\mid\abs{z}<R\}$ since we can use \eqref{2.10} to
define $M(z)$ as a meromorphic function in $\{z\mid 1 < \abs{z}<R\}$ and
\eqref{2.11} says the function has matching boundary values on $\abs{z}=1$.

Moreover, \eqref{2.10} implies
\begin{SL}
\item[(i)] $M(z)$ has a pole at $z_k$ with $1<\abs{z_k} <R$ only if $u(z)$ has a zero
there and the order of the pole is the same as the order of the zero. This is because
$u^\sharp (z)$ is nonvanishing near $z=z_k$ and $M^\sharp(z)$ is regular near $z=z_k$,
since $z_k^{-1} \in\bbD$, and we are supposing no bound states.

\item[(ii)] If $u$ vanishes at $+1$ or $-1$, $M$ has a first-order pole there
(for if $M(z)$ has a pole at $\pm 1$ with real residue, $M^\sharp (z)$ has the
opposite residue, so $M-M^\sharp$ still has a pole).
\end{SL}

Combining this with \eqref{2.4}, we see $u^{(1)}(z)$ is analytic in $\{z\mid\abs{z}<R\}$
for any poles of $M$ are cancelled by zeros of $u$. Moreover, $u^{(1)}$ is nonvanishing
in $\overline{\bbD}$ for $u$ is nonvanishing on $\bbD\backslash \{-1,1\}$ and $M(z)/z$ is
nonvanishing since $\Ima M>0$ on $\bbD\cap\bbC_+$, so $\Real M>0$ for $z\in (0,1)$ and
$\Real M <0$ for $z\in (-1,0)$. It follows $u^{(1)}$ is nonvanishing on
$\overline{\bbD}\backslash \{-1,1\}$. And (ii) above shows that even if $u(\pm 1)$ is
zero, $M$ has a compensating pole.
\end{proof}

We will also need

\begin{theorem} \lb{T2.2} If the Jost function of Jacobi data $\{a_n,
b_n\}_{n=1}^\infty$ has finitely many zeros in $\bbD$ and the only zeros
on $\partial\bbD$ are at $\pm 1$ and those are simple, then
\begin{equation} \lb{2.8a}
\abs{a_n-1} + \abs{b_n} \to 0
\end{equation}
and
\begin{equation} \lb{2.8b}
M^{(n)}(z)\to z
\end{equation}
uniformly on compacts of $\bbD$. In particular, for each $\rho <1$,
\begin{equation} \lb{2.8c}
\sup_{\abs{z}\leq \rho}\, \biggl| \f{M^{(n)}(z)}{z}\biggr| \to 1
\end{equation}
\end{theorem}

\begin{proof} Since the weight of the spectral measure is given by \eqref{1.12x}
and \eqref{1.11}, the Szeg\H{o} condition holds and so does the quasi-Szeg\H{o}
condition of \cite{KS}. This plus finite spectrum show $\sum_{n=1}^\infty
\abs{a_n -1}^2 + \abs{b_n}^2 <\infty$ by the work of Killip-Simon \cite{KS}.
Thus \eqref{2.8a} holds.

That implies the corresponding Jacobi matrix $J^{(n)}$ converges in norm
to $J_0$ so the resolvents converge, which implies \eqref{2.8b}. \eqref{2.8c}
is a consequence of $M^{(n)}(z)/z\to 1$ uniformly.
\end{proof}

We now combine \eqref{2.4} and \eqref{2.10} to write the critical update
equation:
\begin{equation} \lb{2.9a}
u^{(n+1)} (z) = a_{n+1} (1-z^{-2}) (u^{(n)\sharp}(z))^{-1}
+ a_{n+1} z^{-2} u^{(n)}(z) N_n^\sharp(z)
\end{equation}
where
\begin{equation} \lb{2.10a}
N_n(z) = \f{M^{(n)}(z)}{z}
\end{equation}
so
\begin{equation} \lb{2.11a}
N_n^\sharp (z) = zM^{(n)\sharp}(z)
\end{equation}

\eqref{2.9a} looks complicated because of the $u^{(n)\sharp}$ term. But consider
expanding all functions in a Laurent series near $\{z\mid\abs{z}=R_1\}$
for $1<R_1 <R$. $u^{(n+1)}(z)$ only has nonnegative powers and $u^{(n)\sharp}$,
and thus $(1-z^{-2})(u^{(n)\sharp})^{-1}$ only has nonpositive powers. Thus
the first term in \eqref{2.9a} compensates for the negative powers in the
second term, and that is its only purpose. If we project onto positive
powers, it disappears!

We thus define $P_+$ to be the projection in $L^2 (R_1 \partial\bbD, \f{d\theta}{2\pi})$
onto $\{e^{in\theta}\}_{n=1}^\infty$. Applying $P_+$ to \eqref{2.9a}, we find
\begin{equation} \lb{2.12}
u^{(n+1)}(R_1 e^{i\theta}) - u^{(n+1)}(0) = a_{n+1} P_+ [(R_1 e^{i\theta})^{-2}
[u^{(n)}(R_1 e^{i\theta})- u^{(n)}(0)] N_n^\sharp (R_1 e^{i\theta})]
\end{equation}
where we used the fact that $z^{-2} u^{(n)}(0) N_n^\sharp(z)$ has only
negative Laurent coefficients.

Define
\begin{equation} \lb{2.13}
|||g|||_{R_1} = \biggl(\int \abs{g(R_1 e^{i\theta})-g_1 (0)}^2 \,
\f{d\theta}{2\pi} \biggr)^{1/2}
\end{equation}
for functions analytic in a neighborhood of $\{z\mid\abs{z}\leq R_1\}$. Since
$P_+$ is a projection in $L^2$, we obtain
\begin{equation} \lb{2.14}
|||u^{(n+1)}|||_{R_1} \leq a_{n+1} R_1^{-2} \|N_n^\sharp (R_1 e^{i\theta})\|_\infty
|||u^{(n)}|||_{R_1}
\end{equation}

\begin{proof}[Proof of Theorem~\ref{T1.3A}] It is easy to see that $u$ is a
polynomial of exact degree $k$ if and only if $u$ is entire and
\begin{equation} \lb{2.16}
\lim_{R\to\infty} \, \f{|||u|||_R}{R^k} \in (0,\infty)
\end{equation}
Thus, by \eqref{2.14}, if $u$ is a polynomial of  degree $\ell$ and $n>\ell/2$,
then $|||u^{(n)}|||_R =0$, that is, $u^{(n)}$ is a constant. But then the weight in
$M$ is the free one, that is, $a_{j+n}\equiv 1$, $b_{j+n}\equiv 0$ for $j\geq 0$.
Of course, $L$ is a polynomial if and only if $u$ is.
\end{proof}

{\it Remark.} This proof and the direction in the appendix allow us to relate
the degree of the polynomial $u$ to the support of $J-J_0$.

\smallskip
\eqref{2.14} also implies

\begin{proposition} \lb{P2.3} For $1<R_1 <R$, we have
\begin{equation} \lb{2.15}
\limsup \, |||u^{(n)}|||_{R_1}^{1/n} \leq R_1^{-2}
\end{equation}
\end{proposition}

\begin{proof} Note first that
\begin{align*}
\sup_\theta\, \abs{N_n^\sharp (R_1 e^{i\theta})} &= \sup_\theta
\abs{N_n (R_1^{-1} e^{i\theta})} \\
&\leq \sup_{\abs{z}\leq R_1^{-1}}\, \biggl|\f{M^{(n)}(z)}{z}\biggr|
\end{align*}
which goes to $1$ by \eqref{2.8c}. Since $a_n\to 1$ by \eqref{2.8a},
\[
\lim_{n\to\infty}\, \biggl( \, \prod_{j=0}^{n-1} a_{j+1} \|N_j^\sharp
(R_j e^{i\theta})\|_\infty \biggr)^{1/n} =1
\]
so \eqref{2.14} implies \eqref{2.15}.
\end{proof}

\begin{proof}[Proof of Theorem~\ref{T1.4}] If $f(z) =\sum_{n=1}^\infty
A_n z^n$, then
\[
|||f|||_{R_1}^2 = \sum_{n=1}^\infty\, \abs{A_n}^2 R_1^{2n}
\]
is monotone in $R_1$, so \eqref{2.15} implies, by taking $R_1\to R$,
$\limsup |||u^{(n)}|||_{1+\veps}^{1/n} \leq R^{-2}$. Since the Cauchy
integral formula shows
\[
\sup_{\abs{z}\leq 1}\, \abs{f(z) -f(0)} \leq |||f|||_{1+\veps}
\]
we see that, for every $\delta >0$,
\begin{equation} \lb{2.17}
\sup_{\abs{z}\leq 1}\, \abs{u^{(n)}(z) -u^{(n)}(0)}\, \abs{R-\delta}^{-2n}
\to 0
\end{equation}
which in turn, using $u^{(n)}(0)\to 1$ (by \eqref{2.8b}), implies
\[
1=2\int_0^\pi \f{\sin^2\theta}{\abs{u^{(n)}(e^{i\theta})}^2}\, d\theta
= \f{1}{\abs{u^{(n)}(0)}^2} + O(\abs{R-\delta}^{-2n})
\]
which implies
\begin{equation} \lb{2.18}
u^{(n)}(0) =1+O(\abs{R-\delta}^{-2n})
\end{equation}

Thus the difference between the free weight and the weight for $f^{(n)}$ is
$O(\abs{R-\delta}^{-2n})$, so
\begin{equation} \lb{2.19}
\limsup \biggl(\, \sup_{\abs{z}\leq\f12}\, \biggl|\f{M^{(n)}(z)}{z} -1\biggr|
\biggr)^{1/n} \leq R^{-2}
\end{equation}
By \eqref{2.8}
\begin{align*}
\limsup \abs{b_n}^{1/n} &\leq R^{-2} \\
\limsup \abs{(a_{n+1}^2-1) + \tfrac12\, b_{n-1}^2}^{1/n} & \leq R^{-2}
\end{align*}
which implies \eqref{1.18}.
\end{proof}

That completes what we want to say about OPRL with no bound states. As an aside,
we show how the ideas of this section provide an alternate to the hard (i.e.,
inverse spectral) side of the Nevai-Totik theorem, Theorem~\ref{T1.1}. Their
proof is shorter but relies on a magic formula (see (2.4.36) of \cite{OPUC1})
\[
d\mu_\s = 0\Rightarrow \alpha_n = -\kappa_\infty \int \ol{\Phi_{n+1}(e^{i\theta})}\,
D(e^{i\theta})^{-1}\, d\mu(\theta)
\]
Our proof will exploit or develop the relative Szeg\H{o} function, $\delta_0 D$,
of Section~2.9 of \cite{OPUC1}. Our goal is to prove

\begin{theorem}\lb{T2.4} Let $d\mu$ be a measure on $\partial\bbD$ with
$d\mu_\s =0$ and so that the Szeg\H{o} condition holds. Suppose $D(z)^{-1}$
has an analytic continuation to $\{z\mid\abs{z}<R\}$ for some $R>1$. Then
\begin{equation} \lb{2.20}
\limsup_{n\to\infty}\, \abs{\alpha_n}^{1/n} \leq R^{-1}
\end{equation}
\end{theorem}

So we suppose the Szeg\H{o} condition holds, which is equivalent to
\begin{equation} \lb{2.21}
\sum_{n=0}^\infty \, \abs{\alpha_n}^2 <\infty
\end{equation}
Let $d\mu_n$ be the measure with Verblunsky coefficients $\{\alpha_{k+n}\}_{k=0}^\infty$
and $D^{(n)}$ its Szeg\H{o} function. $F^{(n)}$ and $f^{(n)}$ are defined by
\begin{align}
F^{(n)}(z) &= \int \f{e^{i\theta}+z}{e^{i\theta}-z}\, d\mu^{(n)}(\theta) \lb{2.22} \\
F^{(n)}(z) &= \f{1+z f^{(n)}(z)}{1-zf^{(n)}(z)}  \lb{2.23}
\end{align}
Geronimus' theorem (see \cite{OPUC1}) says that the relation between the $f$'s
is given by the Szeg\H{o} algorithm,
\begin{equation} \lb{2.24x}
f^{(n)}(z) \equiv \f{\alpha_n + zf^{(n+1)}(z)}{1+\bar\alpha_n zf^{(n+1)}(z)}
\end{equation}
and the equivalent
\begin{equation} \lb{2.25x}
zf^{(n+1)}(z) = \f{f^{(n)}(z)-\alpha_n}{1-\bar\alpha_n f^{(n)}(z)}
\end{equation}

In Section~2.9 of \cite{OPUC1}, the relative Szeg\H{o} function is defined by
($\rho_n = (1-\abs{\alpha_n}^2)^{1/2}$)
\begin{equation} \lb{2.24}
(\delta_nD) (z) = \f{1-\bar\alpha_n f^{(n)}(z)}{\rho_n} \, \,
\f{1-zf^{(n+1)}(z)}{1-zf^{(n)}(z)}
\end{equation}
and it is proven that
\begin{equation} \lb{2.25}
(\delta_n D)(z) = \f{D^{(n)}(z)}{D^{(n+1)}(z)}
\end{equation}
which we write as
\begin{equation} \lb{2.26}
D^{(n+1)}(z)^{-1} = D^{(n)}(z)^{-1} (\delta_n D)(z)
\end{equation}

It will be useful to rewrite \eqref{2.24} using \eqref{2.25x} to get
\begin{equation} \lb{2.27}
(\delta_n D)(z) = \f{1-\bar\alpha_n f^{(n)}(z) - f^{(n)}(z) + \alpha_n}
{\rho_n (1-zf^{(n)}(z))}
\end{equation}

Using
\begin{equation} \lb{2.27a}
f^{(n)}(z) = \f{1}{z}\, \f{F^{(n)}(z)-1}{F^{(n)}(z)+1}
\end{equation}
one finds
\begin{equation} \lb{2.28b}
(\delta_n D)(z) = \tfrac12\, z^{-1} M^{(n)}(z)
\end{equation}
where
\begin{equation} \lb{2.28c}
M^{(n)}(z) = z(1+\alpha_n) (F^{(n)}(z) +1) - (1+\bar\alpha_n)
(F^{(n)}(z) -1)
\end{equation}
Interestingly enough, $M^{(n)}(z)$ for $n=0$ appears in the theory of minimal
Carath\'eodory functions on the hyperelliptic Riemann surfaces that occur in the analysis
of OPUC with periodic Verblunsky coefficients (see (11.7.76) in \cite{OPUC2}); a related
function appears in Geronimo-Johnson \cite{GJo1}. While \cite{OPUC1,OPUC2} introduced
both $\delta_0D$ and $M^{(0)}(z)$, its author appears not to have realized the relation
\eqref{2.28b}. $\delta_nD$ is nonsingular at $z=0$ since $M^{(n)}(z) = 0$ at $z=0$ (since
$F^{(n)}(0)=1$). We note that where Section~11.7 of \cite{OPUC2} uses $M(z)$, it could
use $(\delta_0 D)(z)$. The difference is the $0$ at $0_+$ is moved to $\infty_-$ and the
pole at $\infty_+$ to $0_-$. We note that the relation between $M$ and $\delta_0D$ is
hinted at in \eqref{2.24}. In gaps in $\supp (d\mu)$ in $\partial\bbD$, $\delta_0D$ has
poles at zeros of $1-zf$ and zeros at zeros of $1-zf_1$. By \eqref{2.23}, $\delta_0 D$
has poles at poles of $F$ and zeros at poles of $F^{(1)}$, which is the critical property
that $M$ needs in the analysis of Section~11.7 of \cite{OPUC2}.

We will also need the analytic continuation of
\begin{equation} \lb{2.28}
\Real F^{(n)} (e^{i\theta}) = \abs{D^{(n)}(e^{i\theta})}^2
\end{equation}
namely,
\begin{equation} \lb{2.29}
F^{(n)} + (F^{(n)})^\sharp = 2D^{(n)} (D^{(n)})^\sharp
\end{equation}
where $\sharp$ is given by \eqref{2.9}.

\begin{theorem}\lb{T2.5} Let $R>1$. If $D^{-1}$ is analytic in $\{z\mid
\abs{z}<R\}$, then $F$ is meromorphic there with singularities precisely at the zeros of
$D^{-1}$. The order of any pole of $F$ is precisely the same as the order of the zero of
$D^{-1}$. $(\delta_0 D)(z)$ is meromorphic in the region with poles precisely at the
poles of $F$ with order no greater than those of $F$\!. $(D^{(1)})^{-1}$ is analytic in
$\{z\mid\abs{z} <R\}$ and thus, by induction, so is each $(D^{(n)})^{-1}$.
\end{theorem}

\begin{proof} By \eqref{2.29}, we can use
\begin{equation} \lb{2.30}
\ti F = -F^\sharp + 2DD^\sharp
\end{equation}
to define a function meromorphic in $\{z\mid R^{-1} < z < R\}$ with poles
in $\{z\mid 1<z<R\}$ precisely at the zeros of $D^{-1}$. Since $\Real F
(e^{i\theta}) = w(e^{i\theta}) = \abs{D(e^{i\theta})}^2$,
\[
\ti F(e^{i\theta}) = - \ol{F(e^{i\theta})} + 2\Real F(e^{i\theta})
= F(e^{i\theta})
\]
$\ti F$ and $F$ agree there and so $\ti F$ extends $F$ to a meromorphic
function in the required region.

By \eqref{2.28b} and \eqref{2.28c}, $(\delta_0 D)(z)$ is meromorphic in
the same region with poles of order no greater than those of $F$\!. Thus
$D^{-1} (\delta_0D)$ is analytic in $\{z\mid \abs{z}<R\}$.
\end{proof}

Now define
\begin{align}
A(z) &= \tfrac12 \left[ z(1+\alpha_0) (1-F^\sharp(z)) + (1+\bar\alpha_0)
(1+F^\sharp (z)) \right] \lb{2.31} \\
B(z) &= (1+\alpha_0) D^\sharp - z^{-1} (1+\bar\alpha_0) D^\sharp \lb{2.32}
\end{align}
so, by \eqref{2.28b}, \eqref{2.28c}, and \eqref{2.29},
\begin{equation} \lb{2.33}
\delta_0 D = z^{-1} A+BD
\end{equation}
and thus
\begin{equation} \lb{2.34x}
[D^{(1)}]^{-1} = z^{-1} A (D^{-1}-D(0)^{-1}) + B+z^{-1} AD(0)
\end{equation}

\begin{proof}[Proof of Theorem~\ref{T2.4}] For functions continuous on $\{z\mid
\abs{z}=R_1\}$, define $P_+$ to be the projection onto positive Fourier terms
and $|||g||| = \|P_+ g\|_{L^2}$. Then \eqref{2.34x} implies, for any $R_1 <R$,
\begin{equation} \lb{2.34}
|||[D^{(1)}]^{-1}|||_{R_1} \leq R_1^{-1}\, \sup_{\abs{z}=R_1}\, \abs{A(z)}\,\,
|||D^{-1}|||_{R_1}
\end{equation}

By induction,
\begin{equation} \lb{2.35}
|||[D^{(n)}]^{-1}|||_{R_1} \leq R_1^{-n} a_0 \dots a_{n-1} |||D^{-1}|||_{R_1}
\end{equation}
where
\begin{equation} \lb{2.36}
a_j =\sup_{\abs{z}=R_1}\, \abs{A^{(j)}(z)} \leq \tfrac12\, (1+\abs{\alpha_j})\,
\sup_{\abs{z}\leq R_1^{-1}}\, (\abs{(1-F^{(j)}(z))z^{-1}} + \abs{1+F^{(j)}(z)})
\end{equation}

Since $\sum_{n=0}^\infty \abs{\alpha_n}^2 <\infty$, we have $\sup_{\abs{z}\leq R_1^{-1}}
\, \abs{f^{(j)}(z)}\to 0$, which implies, by \eqref{2.27a}, that
\[
\sup_{\abs{z}\leq R_1^{-1}}\, \abs{(1-F^{(j)}(z))z^{-1}}\to 0
\]
and thus, by \eqref{2.36},
\[
\limsup_{j\to\infty}\, \abs{a_j} \leq 1
\]
so \eqref{2.35} implies
\[
\limsup_{n\to\infty}\, ||| [D^{(n)}]^{-1}|||_{R_1}^{1/n} \leq R_1^{-1}
\]

By analyticity, this implies (taking $R_1\to R$)
\[
\limsup_{n\to\infty}\, \|w_n -1\|_\infty^{1/n} \leq R^{-1}
\]
which implies
\[
\lim_{n\to\infty}\, \abs{\alpha_n}^{1/n} \leq R^{-1}
\]
since
\[
\alpha_n = \int e^{-i\theta} (w_n (\theta) -1) \, \f{d\theta}{2\pi}
\qedhere
\]
\end{proof}

\section{The General Case} \lb{s3}

In this section, we will prove Theorems~\ref{T1.5} and \ref{T1.6}. A key piece of
the proofs is the following well-known result:

\begin{theorem}\lb{T3.1} If $\{a_n, b_n\}_{n=1}^\infty$ are a set of Jacobi parameters
with $\spec (J)\backslash [-2,2]$ finite, then there is a $k$ so that $\{a_{n+k},
b_{n+k}\}_{n=1}^\infty$ are a set of Jacobi parameters with $\spec(J_k) \subset [-2,2]$.
\end{theorem}

\begin{proof} By a Sturm oscillation theorem (see \cite{Geneva,CL}), the number of
spectral points in $(2,\infty)$ is the number of sign flips of $\{P_n(z)\}_{n=1}^\infty$
and in $(-\infty, -2)$ of $\{(-1)^n P_n(-2)\}_{n=1}^\infty$. By assumption, these
numbers are finite, so for some $k$, $\{P_{n+k-1}(2)\}_{n=1}^\infty$ and $\{(-1)^n
P_{n+k-1}(-2)\}_{n=1}^\infty$ have fixed signs. By a comparison theorem and the
oscillation theorem again, it follows that $\spec (J_k) \subset [-2,2]$.
\end{proof}

The reason weights have to be canonical is that Theorem~\ref{T2.1} can fail if there are
zeros in $\bbD$. This is because $u^{(1)} = a_1 z^{-1} u M$ and
\begin{equation}  \lb{3.1}
M(z) =M^\sharp(z) + (z-z^{-1}) [u(z) u^\sharp(z)]^{-1}
\end{equation}
If $z_j$ is a (real) zero in $\bbD$ with $\abs{z_j} >R^{-1}$, $M(z)$ may have a
pole at $1/z_j$ due to the $(u^\sharp)^{-1}$ term. $uM$ will then have a pole
(unless $u$ has a zero, which we will see does not help). The way to avoid this
is to arrange for $M^\sharp$ to have a compensating pole, and this will happen
precisely if the weight at $z_j$ is canonical! Here is the detailed result:

\begin{theorem}\lb{T3.2} Let $u$ be a Jost function and be analytic in $\{z\mid
\abs{z}<R\}$ and suppose there are finitely many zeros $\{z_j\}_{j=1}^N$ in $\bbD$. Then
$M$ is meromorphic in $\{z\mid\abs{z}<R\}$ and the only possible poles of $M$ in $\{z\mid
1 < \abs{z} < R\}$ are at points $z_j^{-1}$ where $\abs{z_j} > R^{-1}$. Moreover,
\begin{SL}
\item[{\rm{(i)}}] If $u$ has a zero of order $k\geq 1$ at $z_j^{-1}$,
then $M$ has a pole of order $k+1$ there.
\item[{\rm{(ii)}}] If $u(z_j^{-1})\neq 0$, then $M$ has a pole at $z_j^{-1}$
if and only if the weight at $z_j$ is {\em not} canonical.
\end{SL}
In particular, $u^{(1)} =uM$ is analytic in $\{z\mid\abs{z}<R\}$ if and
only if all weights at those $z_j$ with $\abs{z_j}>R^{-1}$ are canonical.
\end{theorem}

\begin{proof} The zeros at $z_j$ are simple, so if $u$ has a zero of order
$k\geq 1$ at $z_j^{-1}$, then $(u(z) u^\sharp (z))^{-1}$ has a pole of
order $k+1$. Since $M^\sharp(z)$ has only simple poles in $\{z\mid\abs{z}
>1\}$, $M(z)$ has a pole of order $k+1$. This proves (i).

If $u(z_j^{-1})\neq 0$, both $(z-z^{-1}) [u(z)u^\sharp(z)]^{-1}$ and $M^\sharp (z)$ have
simple poles at $z_j^{-1}$. Their residues cancel if and only if \eqref{1.20} holds.
\end{proof}

This gets us one step if the weights are canonical. We get beyond that because
automatically weights after that are canonical!

\begin{theorem}\lb{T3.3} Let $u$ be a Jost function and be analytic in $\{z\mid
\abs{z}<R\}$. Let $u^{(1)}= a_1 z^{-1} uM$ and let $M^{(1)}$ obey \eqref{2.5}.
Then $\ti z_j \in\bbD$ with $\abs{\ti z_j}>R^{-1}$ implies that $u^{(1)}
(\ti z_j^{-1})\neq 0$ and the weight of $M^{(1)}$ is canonical.
\end{theorem}

\begin{proof} Zeros of $u$ in $\bbD$ are cancelled by poles of $M$\!, so
$u^{(1)}$ has zeros precisely at points $\ti z_j$ where $M(\ti z_j)=0$. Since $[M(\ti
z_j) -M^\sharp(\ti z_j)] u^\sharp (\ti z_j) u (\ti z_j) = \ti z_j -\ti z_j^{-1} \neq 0$,
$M^\sharp (\ti z_j) \neq 0$. By \eqref{2.5}, $M(\ti z_j^{-1}) \neq 0$ implies
$M^{(1)}(z)$ is regular at $\ti z_j^{-1}$. By a small calculation, $u^{(1)}, M^{(1)}$
obey \eqref{2.10},  so the weight must be canonical.
\end{proof}

{\it Remark.} There is a potentially puzzling feature of Theorem~\ref{T3.3}.
If stripping a Jacobi parameter pair cannot produce noncanonical weights,
how can they occur? After all, we can add a parameter pair before $J$ and
then remove it. The resolution is that adding a parameter pair also shrinks
the region of analyticity if $J$ has a noncanonical weight. Essentially,
noncanonical weights produce poles in a meromorphic $u^{(n)}(z)$ with
residues which are a ``resonance eigenfunction." Lack of a canonical
weight in $u^{(n)}$ is a sign that this resonance eigenfunction function
vanishes at $n$, but then it will not at $n+1$ or $n-1$.

\begin{proof}[Proof of Theorem~\ref{T1.5}] If some weight is not canonical,
$u^{(1)}$ is not entire, and so the Jacobi parameters cannot have finite support. If all
weights are positive, one can normalize $u$ so \eqref{1.13} holds with the canonical
weights and obtain a positive weight since the necessary canonical weights are all
positive. With canonical weights, $uM$ is entire and by \eqref{3.1} and $M(z) =O(z)$ at
$z=0$, we see that $M(z) =O(1/z)$ at $z=\infty$ so long as $u(z) =O(z^\ell)$ with $\ell
>2$. Thus $u^{(1)}$ is a polynomial of degree at most $1$ less than $u$. Iterating and
using Theorem~\ref{T3.1}, we eventually get a polynomial of $u^{(k)}(z)$ with no zeros in
$\bbD$ and so, by Theorem~\ref{T1.3A}, a finite-range set of $\{a_n, b_n\}$.
\end{proof}

\begin{proof}[Proof of Theorem~\ref{T1.6}] If weights are canonical, we can
iterate to a $u^{(k)}$ nonvanishing in $\bbD$, and use Theorem~\ref{T1.4}.
If some weight is not canonical, $u^{(1)}$ is not analytic in $\{z\mid\abs{z}
<R\}$, so \eqref{1.18} fails on account of Theorem~\ref{TA.1.3}.
\end{proof}

\appendix
\section{The Geronimo-Case Equations} \lb{App}
\renewcommand{\theequation}{A.\arabic{equation}}
\renewcommand{\thetheorem}{A.\arabic{theorem}}
\setcounter{theorem}{0}
\setcounter{equation}{0}

In this appendix, we provide a proof of basic facts about the Jost
functions in the case of regular bounds on the Jacobi coefficients.
We do this primarily because we want to make propaganda for a lovely
set of equations of Geronimo-Case \cite{GC80} which have not yet gotten
the attention they deserve. It is also useful to keep this paper
self-contained. More well-known approaches to Jost functions involve
finding the Jost solution by using variation of parameters about the
free solutions (see, e.g., Teschl \cite{Teschl}) or, going back to
Jost-Pais \cite{JP}, as perturbation determinants (see, e.g.,
Killip-Simon \cite{KS}).

There is some overlap in our presentation and that of Geronimo-Case
\cite{GC80} and Geronimo \cite{Ge94}, but we feel it might be useful
to present the detailed estimates concisely. And we wish to emphasize
the a priori derivation of the GC equations (rather than presenting them
and proving they have asymptotic properties) by identifying what their
$\psi_n$ ($=z^{-n} g_n(z)$ for our $g_n$ below) is. In any event, we
have put this material in an appendix since we regard it as a review.

Given a Jacobi matrix, as in \eqref{1.1a}, we define $a_0=1$ and look
at solutions $(f_n)_{n=0}^\infty$ of
\begin{equation} \lb{A.1.1}
a_n f_{n+1} + (b_n -E) f_n + a_{n-1} f_{n-1} =0 \qquad n=1,2,\dots
\end{equation}

Of course, one solution of this is
\begin{equation} \lb{A.1.2}
f_n = p_{n-1} (E,J)
\end{equation}
with $p_n$ the orthonormal polynomials associated to $J$.

If $f_n$ and $k_n$ are two sequences, we define their Wronskian (we will
often drop the $J$),
\begin{equation} \lb{A.1.3}
W(f,k; J)(n) \equiv a_n (f_{n+1} k_n - f_n k_{n+1})
\end{equation}
If $f,k$ both solve \eqref{A.1.1}, the Wronskian is constant.

Given a Jacobi matrix $J$, define $\ti J_\ell$ by setting $\ti a_{\ell +1}
=\ti a_{\ell+2} = \cdots =1$ and $\ti b_{\ell+1} = \ti b_{\ell+2} =\cdots
= 0$ (in \cite{KS} and \cite{SZ}, $J_\ell$ is defined also setting $a_\ell
=1$; it is different from $\ti J_\ell$), that is,
\begin{equation} \lb{A.1.4}
\ti J_\ell = \begin{pmatrix}
b_1 & a_1 & 0 & {} & {} & {} & {} & {} & {} \\
a_1 & b_2 & a_2 & {} & {} & {} & {} & {} & {} \\
0 & a_2 & \ddots & \ddots & {} & {} & {} & {} & {} & {} \\
{} & {} & \ddots & \ddots & \ddots \\
{} & {} & {} & a_{\ell-1} & b_\ell & a_\ell & 0 \\
{} & {} & {} & 0 & a_\ell & 0 & 1 \\
{} & {} & {} & 0 & 0 & 1 & 0 \\
{} & {} & {} & {} & {} & \ddots & \ddots & \ddots \\
{} & {} & {} & {} & {} & {} & \ddots & \ddots & \ddots
\end{pmatrix}
\end{equation}

$u_n (z; \ti J_\ell)$ solves \eqref{A.1.1} for $\ti J_\ell$ with $u_n (z;
\ti J_\ell)=z^n$ if $n\geq \ell+1$ where $z+z^{-1} =E$. $u_n$ is called
the {\it Jost solution}. The {\it Jost function\/} is
\begin{align}
u(z;\ti J_\ell) &= W(p_{\boldsymbol{\cdot}-1} (z+z^{-1}; \ti J_\ell),
u_{\boldsymbol{\cdot}} (z, \ti J_\ell)) \lb{A.1.5} \\
&= u_0 (z;\ti J_\ell) \notag
\end{align}
by taking the Wronskian at $n=1$. Since $u_\ell (z;\ti J_\ell) = a_\ell^{-1}
z^\ell$, taking the Wronskian at $n=\ell$, we find
\begin{equation} \lb{A.1.6}
u(z; \ti J_\ell) = a_\ell (p_\ell (z+z^{-1}; \ti J_\ell) \, \f{1}{a_\ell}\,
z^\ell - p_{\ell-1} (z+z; \ti J_\ell) z^{\ell+1})
\end{equation}

Since the parameters of $\ti J_\ell$ and $J$ agree for $a_m, b_m$, $m=1,2,\dots,
\ell$, we have $p_m (z,J) = p_m (z,J_\ell)$ for $m\leq \ell$. Thus we define
\begin{equation} \lb{A.1.7}
g_n(z) = z^n \biggl( p_n \biggl( z+\f{1}{z}\biggr) - a_n zp_{n-1} \biggl( z+
\f{1}{z}\biggr)\biggr)
\end{equation}
which is a polynomial in $z$ of degree at most $z^{2n}$. $g_n(z)$ is the
Jost function for $\ti J_n$, that is,
\begin{equation} \lb{A.1.8}
g_n(z) = u(z;\ti J_n)
\end{equation}
It is also natural to define
\begin{equation} \lb{A.1.9}
c_n (z) =z^n p_n \biggl( z+\f{1}{z}\biggr)
\end{equation}
a polynomial in $z$ of exact degree $2n$.

Clearly, \eqref{A.1.7} becomes
\begin{equation} \lb{A.1.10}
g_{n+1}(z) = c_{n+1}(z) - a_{n+1} z^2 c_n (z)
\end{equation}
The fundamental equation
\[
a_{n+1} p_{n+1} \biggl( z+\f{1}{z}\biggr) = \biggl( z +\f{1}{z} -
b_{n+1}\biggr) p_n \biggl(z + \f{1}{z}\biggr) - a_n p_{n-1}
\biggl( z + \f{1}{z}\biggr)
\]
multiplied by $z^{n+1}$ becomes (using \eqref{A.1.10} for $n\to n-1$)
\begin{equation} \lb{A.1.11}
a_{n+1} c_{n+1} = (z^2 -b_{n+1} z) c_n + g_n
\end{equation}
Multiplying \eqref{A.1.10} by $a_{n+1}$ and using \eqref{A.1.11},
we find
\begin{equation} \lb{A.1.12}
a_{n+1} g_{n+1} = [(1-a_{n+1}^2) z^2 - b_{n+1} z]c_n + g_n
\end{equation}

The last two equations, which we call the Geronimo-Case equations, are
the fundamental recursion relations with initial conditions
\begin{equation} \lb{A.1.13}
g_0 (z) = c_0 (z) =1
\end{equation}
(One checks $g_0 =1$ by looking at \eqref{A.1.11} using $c_0 =1$ and $c_1 = \f{z}{a_1}
(z+\f{1}{z} -b_1)$.) It is natural to define the Geronimo-Case update matrix
\begin{equation} \lb{A.1.14}
U_n(z) = \begin{pmatrix}
z^2 - b_n z & 1 \\ (1-a_n^2) z^2 - b_n z & 1 \end{pmatrix}
\end{equation}
so \eqref{A.1.11}/\eqref{A.1.12} become
\begin{equation} \lb{A.1.15}
\binom{c_{n+1}}{g_{n+1}} = \f{1}{a_{n+1}} \, U_{n+1} \binom{c_n}{g_n}
\end{equation}
or, with
\begin{gather}
T_n = U_n U_{n-1} \dots U_1  \lb{A.1.16} \\
\binom{c_n(z)}{g_n(z)} = (a_1 \dots a_n)^{-1} T_n (z) \binom{1}{1} \lb{A.1.17}
\end{gather}

If we use unnormalized functions
\[
C_n(z) = a_1 \dots a_n c_n (z) \qquad G_n(z) = a_1 \dots a_n g_n (z)
\]
we have
\begin{equation} \lb{A.1.18}
\binom{C_n}{G_n} = T_n \binom{1}{1}
\end{equation}
Since we have the normalization conditions $C_n(z) =z^{2n} +$ lower order, $G_n(0)=1$,
this has some similarity to the equations for orthogonal polynomials on the unit circle
\cite{OPUC1,OPUC2,Szb,GBk, GBk1}, a motivation for Geronimo-Case \cite{GC80}, but there
are numerous differences; for example, $C_n^*(z) = G_n(z)$ for the reversal operations of
the Szeg\H{o} recursion but not here.

One consequence of this is immediate since $u=g_n$ if $a_\ell =1$, $b_\ell =0$ for
$\ell\geq n+1$:

\begin{theorem}\lb{TA.1.1} If $J-J_0$ is finite range, then $u(z;J)$ is a
polynomial. If $a_\ell =1$, $b_\ell =0$ for $\ell\geq n+1$ and $a_n\neq 1$,
then $\deg (u) =2n$; if $a_n =1$ but $b_n\neq 0$, $\deg (u) =2n-1$.
\end{theorem}

{\it Remark.} Indeed, in the case $a_n\neq 1$, the proof shows that
\begin{equation} \lb{A.1.19}
u(z;J) = \f{(1-a_n^2)}{a_1\dots a_n}\, z^{2n} + \text{ lower order}
\end{equation}
and if $a_n=1$ but $b_n\neq 0$,
\begin{equation} \lb{A.1.20}
u(z;J) = -\f{b_n}{a_1 \dots a_n} \, z^{2n-1} + \text{ lower order}
\end{equation}

\begin{proof} We have, if $a_\ell =1$, $b_\ell =0$ for $\ell \geq n+1$, that
\[
u(z;J) = a_n^{-1} [(1-a_n^2) z^2 - b_n z] c_{n-1}(z) + g_{n-1}(z)
\]
Since $\deg (g_{n-1})\leq 2n-2$ and $c_{n-1}(z) = (a_1 \dots a_{n-1})^{-1}
z^{2n-2} +$ lower order, the result is immediate.
\end{proof}

Notice if $b_n=0$, $a_n=1$, then
\begin{equation} \lb{A.1.21x}
U_n = U^{(0)} (z) \equiv \begin{pmatrix} z^2 & 1 \\ 0 & 1 \end{pmatrix}
\end{equation}

We will make one of three successively stronger hypotheses on the
Jacobi coefficients:
\begin{align}
\sum_{n=1}^\infty \, [\abs{b_n} + \abs{a_n^2 -1}] < \infty \lb{A.1.21} \\
\sum_{n=1}^\infty n [\abs{b_n}+\abs{a_n^2 -1}] <\infty \lb{A.1.22} \\
\abs{b_n} + \abs{a_n^2 -1} \leq CR^{-2n} \lb{A.1.23}
\end{align}
for some $R>1$. $a_n^2 -1$ will enter in estimates, but since $\abs{a_n-1}
\leq \abs{a_n^2-1} \leq (1+\sup_n \abs{a_n}) (\abs{a_n-1})$ in all these
estimates, $\abs{a_n-1}$ can replace $\abs{a_n^2 -1}$ with no change.

In all cases, $\sum \abs{a_n-1}<\infty$, so $\prod_{j=1}^n a_j$ is uniformly bounded
above and below, and hence $c$ (resp., $g$) and $C$ (resp., $G$) are comparable.

We will first prove bounds and use them to control convergence:

\begin{theorem} \lb{TA.1.2}
\begin{SL}
\item[{\rm{(i)}}] Let \eqref{A.1.21} hold. Then for each $z\in\overline{\bbD} \backslash
\{\pm 1\}$,
\begin{equation} \lb{A.1.24}
\sup_n \, [\abs{G_n(z)} + \abs{C_n(z)}] \equiv A_0(z) < \infty
\end{equation}
where $A_0(z)$ is bounded uniformly on compact subsets of $\overline{\bbD}
\backslash\{\pm 1\}$. \item[{\rm{(ii)}}] Let \eqref{A.1.22} hold. Then for some constant
$A_1$,
\begin{align}
\sup_{n,z\in\overline{\bbD}}\, \abs{G_n(z)} &\leq A_1 \lb{A.1.25} \\
\sup_{n,z\in\overline{\bbD}}\, \f{\abs{C_n (z)}}{1+n} &\leq A_1 \lb{A.1.26}
\end{align}
\item[{\rm{(iii)}}] Let \eqref{A.1.23} hold and let $R>1$. Then there is some constant
$A_2$ such that for all $z$ with $\abs{z}<R$,
\begin{equation} \lb{A.1.26a}
\abs{G_n(z)} + \abs{C_n (z)} \leq A_2 (1+n) [\max(1,|z|)]^{2n}
\end{equation}
\end{SL}
\end{theorem}

\begin{proof}
(i) If $z\neq \pm 1$, $U^{(0)}(z)$ can be diagonalized; explicitly,
\begin{equation} \lb{A.1.29}
L(z) U^{(0)}(z) L(z)^{-1} = \begin{pmatrix}  z^2 & 0 \\ 0 & 1 \end{pmatrix}
\end{equation}
where
\begin{equation} \lb{A.1.30}
L(z) = \begin{pmatrix} 1 & -\f{1}{1-z^2} \\ 0 & 1 \end{pmatrix}
\end{equation}
Thus for $z\neq\pm 1$,
\[
L(z) T_n(z) L(z)^{-1} = (K_0(z) + B_n(z)) \dots (K_0(z) + B(z))
\]
with
\[
K_0(z) = \begin{pmatrix} z^2 & 0 \\ 0 & 1 \end{pmatrix} \qquad
B_n(z) = L(z) \begin{pmatrix} -b_n z & 0 \\ -b_n z+(1-a_n^2) z^2 & 0
\end{pmatrix} L(z)^{-1}
\]

If $\abs{z}\leq 1$, $z\neq\pm 1$,
\[
\|K_0(z) + B_n(z)\| \leq 1 + \|L(z)\|\, \|L(z)^{-1} \| \, [\abs{1-a_n^2}
+ \abs{b_n}]
\]
so \eqref{A.1.21} implies \eqref{A.1.24} with
\begin{equation} \lb{A.1.30a}
A_0(z) \le \biggl( 1 + \f{1}{\abs{1-z^2}}\biggr)^2 \prod_{j=1}^\infty \, \biggl[
1+\biggl( 1 + \f{1}{\abs{1-z^2}}\biggr)^2 (\abs{1-a_j^2} + \abs{b_j})\biggr]
\end{equation}

\smallskip
(ii) Let us show inductively that for all $\abs{z}\leq 1$,
\begin{align}
\abs{G_n(z)} &\leq \prod_{j=1}^n \, [1+j (\abs{1-a_j^2} + \abs{b_j})] \lb{A.1.27}  \\
\abs{C_n(z)} &\leq (n+1) \prod_{j=1}^n \, [1+j(\abs{1-a_j^2} + \abs{b_j})] \lb{A.1.28}
\end{align}
This implies \eqref{A.1.25}/\eqref{A.1.26} with
\[
A_1 = \prod_{j=1}^\infty\,  [1+j (\abs{1-a_j^2} + \abs{b_j})]
\]

To prove \eqref{A.1.27}/\eqref{A.1.28}, note they hold for $n=0$ (where they say
$\abs{G_0}\leq 1$, $\abs{C_0}\leq 1$). If they hold for $n$, then by \eqref{A.1.12},
\[
\abs{G_{n+1}(z)} \leq [(n+1)(\abs{a_{n+1}^2-1} + \abs{b_{n+1}}) + 1] \prod_{j=1}^n \,
[1+j(\abs{1-a_j^2} + \abs{b_j})]
\]
and by \eqref{A.1.11},
\[
\abs{C_{n+1}(z)} \leq [(n+1)(1+\abs{b_{n+1}}) + 1] \prod_{j=1}^n \, [1+j(\abs{1-a_j^2} +
\abs{b_j})]
\]

\smallskip
(iii) Since $U_n(z)$ is analytic, the maximum principle says we need only prove the
estimate for $\abs{z} =R-\veps$ with $\veps$ small. As in part (i),
\[
\|K_0(z) +B_n(z)\| \leq \abs{z}^2 [1+\|L(z)\|\, \|L(z)\|^{-1} (\abs{b_n}+\abs{1-a_n^2})]
\]
from which \eqref{A.1.26a} follows.
\end{proof}

\begin{theorem}\lb{TA.1.3}
\begin{SL}
\item[{\rm{(i)}}] Let \eqref{A.1.21} hold. Then
\begin{equation} \lb{A.1.31}
u (z;J) = \lim_{n\to\infty} \, g_n(z)
\end{equation}
converges for all $z\in\overline{\bbD}\backslash \{\pm 1\}$ uniformly on compact subsets
of $\overline{\bbD} \backslash \{\pm 1\}$. $u$ is analytic on $\bbD$, continuous on
$\overline{\bbD}\backslash \{\pm 1\}$. Moreover, for $z\in\bbD$, we have Szeg\H{o}
asymptotics
\begin{equation} \lb{A.1.32}
\lim_{n\to\infty} \, c_n(z) = \f{u(z;J)}{1-z^2}
\end{equation}

\item[{\rm{(ii)}}] Let \eqref{A.1.22} hold. Then \eqref{A.1.31} holds for all
$z\in\overline{\bbD}$ converging uniformly there. $u$ is continuous on $\overline{\bbD}$.
\item[{\rm{(iii)}}] Let \eqref{A.1.23} hold. Then \eqref{A.1.31} holds for
$z\in\{z\mid\abs{z}<R\}$ uniformly on compact subsets.
\end{SL}
\end{theorem}

{\it Remark.} In \cite{Jost1}, we proved in general that if \eqref{A.1.31}
holds, then so does \eqref{A.1.32} for $\abs{z}<1$.

\begin{proof} Again, since $\prod_{j=1}^n
a_j$ is uniformly bounded above and below, it suffices to prove convergence of $G_n$ and
\eqref{A.1.32} in the sense that
\begin{equation} \lb{A.1.33}
\lim_{n\to\infty} \, C_n(z) = (1-z^2)^{-1} \lim_{n\to\infty}\, G_n(z)
\end{equation}
By \eqref{A.1.12},
\begin{equation} \lb{A.1.34}
\abs{G_{n+1}(z) -G_n(z)} \leq (\abs{a_{n+1}^2 - 1}\, \abs{z^2} + \abs{b_{n+1}}\, \abs{z})
\abs{C_n(z)}
\end{equation}
Thus in all cases, the previous theorem implies
\[
\sum_{n=0}^\infty \, \abs{G_{n+1}(z) - G_n(z)}<\infty
\]
uniformly on compacts of the appropriate regions.

All that remains is \eqref{A.1.33}. We have for $\abs{z}<1$,
\[
\abs{C_n - G_{n-1} - z^2 C_{n-1}} \leq A_0(z) \abs{b_n}
\]
Iterating we see that
\[
\biggl| C_n-\sum_{j=0}^{n-1} G_{n-j-1} z^{2j}\biggr| \leq \abs{z}^{2n} + A_0(z)
\sum_{j=0}^{n-1} \, \abs{b_{n-j}}\, \abs{z}^{2j}
\]
For $k \le n-1$, we have
\begin{equation}\lb{A.1.35}
\sum_{j=0}^{n-1} \, \abs{b_{n-j}}\, \abs{z}^{2j} \le \sum_{j=0}^{k-1} \, \abs{b_{n-j}} +
\abs{z}^{2k} \sum_{\ell=k}^{n-1} \, \abs{b_{n-\ell}}
\end{equation}
If we let first $n\to\infty$ and then $k\to\infty$, we see that
$$
\sum_{j=0}^{n-1} \, \abs{b_{n-j}}\, \abs{z}^{2j} \to 0
$$
since $\sum_{k=0}^\infty \abs{b_k} <\infty$. Thus
\[
\lim_{n\to\infty} \, \biggl|C_n -\sum_{j=0}^{n-1} G_{n-j-1}
z^{2j} \biggr| =0
\]

But we claim
\[
\sum_{j=0}^{n-1} (G_{n-j-1} -G_\infty) z^{2j}\to 0
\]
by an inequality like \eqref{A.1.35}. Thus
\[
\biggl| C_n -G_\infty \, \f{1-z^{2n}}{1-z^2}\biggr| \to 0
\]
proving \eqref{A.1.33}.
\end{proof}

To understand where \eqref{A.1.33} comes from, the above actually shows
$T_\infty =\lim_{n\to\infty} T_n$ exists if $z\in\bbD$. Since
\[
T_{n+1}=[U^{(0)}(z) + (U_{n+1}(z) - U^{(0)}(z))] T_n
\]
we see that $U^{(0)}(z) T_\infty (z) = T_\infty (z)$ so $U^{(0}(z)^k
T_\infty (z) =T_\infty (z)$. Since $U^{(0)}(z)$ has eigenvalues $z^2$
and $1$ and $\abs{z^2}<1$, we see that
\[
P(z) T_\infty =T_\infty
\]
with $P(z)$ the spectral projection for eigenvalue $1$, that is, $T_\infty
\binom{1}{1}$ is a multiple of the eigenvectors of $U^{(0)}(z)$ with
eigenvalue $1$. This eigenvector is $\binom{1/(1-z^2)}{1}$, explaining
why \eqref{A.1.33} holds.

If $J$ is a Jacobi matrix, $J^{(\ell)}$ is the Jacobi matrix with
$a_n^{(\ell)}=a_{n+\ell}$, $b_n^{(\ell)} =b_{n+\ell}$. If $J$ obeys
any of \eqref{A.1.21}--\eqref{A.1.23}, so do $J^{(\ell)}$ for all
$\ell$, and thus $u(z;J^{(\ell)})$ exists in the approximate region.
We define the Jost solution (we will show below it agrees with our earlier
definition when $J=\ti J_\ell$) by
\begin{equation} \lb{A.1.36}
u_n (z;J) = a_n^{-1} z^n u(z;J^{(n)})
\end{equation}
By the ``appropriate region" we mean
\begin{alignat*}{2}
&\overline{\bbD}\backslash\{\pm 1\} &&\qquad\text{if \eqref{A.1.21} holds} \\
&\overline{\bbD} &&\qquad\text{if \eqref{A.1.22} holds} \\
&\{z\mid \abs{z}<R \} &&\qquad\text{if \eqref{A.1.23} holds}
\end{alignat*}

\begin{theorem}\lb{TA.1.4}
\begin{SL}
\item[{\rm{(i)}}] $u_n (z;J)$ defined in the appropriate region obeys
\begin{equation} \lb{A.1.37}
a_n u_{n+1} (z;J) + (b_n -(z+z^{-1})) u_n (z;J) + a_{n-1} u_{n-1} (z;J) =0
\end{equation}

\item[{\rm{(ii)}}] In the appropriate region,
\begin{equation} \lb{A.1.37a}
\lim_{n\to\infty} \, z^{-n} u_n (z;J) =1
\end{equation}

\item[{\rm{(iii)}}] For $z\in\bbD$,
\begin{equation} \lb{A.1.38}
u(z;J^{(1)})=a_1 z^{-1} u(z;J) M(z;J)
\end{equation}
where for $z\in\bbD$,
\begin{equation} \lb{A.1.39}
M(z) = (\delta_1, (z+z^{-1} -J)^{-1}\delta_1)
\end{equation}

\item[{\rm{(iv)}}] The only zeros of $u(z;J)$ in $\bbD$ are at real points $\beta$ with
$\beta + \beta^{-1}$ a discrete eigenvalue of $J$. Each zero of $u(z;J)$ in $\bbD$ is
simple.

\item[{\rm{(v)}}] The only zeros of $u(z;J)$ in $\partial\bbD$ are
possible ones at $\pm 1$. If \eqref{A.1.22} holds, these zeros are
simple in the sense that
\begin{equation} \lb{A.1.37b}
u(\pm 1, J)=0 \Rightarrow \lim_{\theta\downarrow 0} \,
\theta^{-1} u(e^{i\theta}, J)\neq 0
\end{equation}

\item[{\rm{(vi)}}] In case \eqref{A.1.23} holds and $u(z_0;J)=0$ with
$R^{-1}<\abs{z_0}<1$, then $u(z_0^{-1};J)\neq 0$.

\item[{\rm{(vii)}}] If \eqref{A.1.21} holds, $M(z)$ has a continuation from $\bbD$ to
$\overline{\bbD}\backslash\{\pm 1\}$, which is everywhere finite and nonzero on
$\partial\bbD\backslash\{\pm 1\}$ and
\begin{equation} \lb{A.1.39a}
\abs{u(e^{i\theta})}^2 \Ima M(e^{i\theta})=\sin\theta
\end{equation}
\end{SL}
\end{theorem}

\begin{proof} (i) Since $u(z;\ti J_\ell) =g_\ell(z;J)\to u(z;J)$, it suffices
to prove \eqref{A.1.37} for $J\equiv \ti J_\ell$ for the {\it new} definition
of $u$. Let $\ti u_n(z;\ti J_\ell)$ temporarily denote the old definition,
that is, the solution of \eqref{A.1.37} with $\ti u_n (z;\ti J_\ell)=z^n$
for large $n$. If we prove $u_n (z;\ti J_\ell)=\ti u_n (z; \ti J_\ell)$,
clearly \eqref{A.1.37} holds for the new definition.

Since $J^{(k)}$ shifts by $k$ steps and $z^n =z^{-k} (z^{n+k})$, we have
for all $n\geq 1$ and $k\geq 1$,
\[
\ti u_n (z; [\ti J_\ell]^{(k)}) =z^{-k} \ti u_{n+k} (z; \ti J_\ell)
\]
Because $\ti u_0 (z;[\ti J_\ell]^{(k)})$ is computed using $\ti a_0=1$
but $\ti u_n (z;\ti J_\ell)$ with $\ti a_0 =a_k$, we have
\[
\ti u_0 (z; [\ti J_\ell]^{(k)}) =a_k z^{-k} \ti u_k (z;\ti J_\ell)
\]
Thus
$$
\ti u_k (z;\ti J_\ell) = a_k^{-1} z^k u(z;[\ti J_\ell]^{(k)}) = u(z;\ti J_\ell)
$$

\smallskip
(ii) Clearly, \eqref{A.1.37a} is equivalent to
\[
\lim_{n\to\infty} \, G_\infty(z;J^{(n)})=1
\]
It follows from \eqref{A.1.34} that
\[
\abs{G_\infty(z;J^{(n)}) -1} \leq \sum_{j=n+1}^\infty \bigl[ \abs{1-a_j^2} \abs{z}^2 +
\abs{b_j} \abs{z} \bigr] \abs{C_{j-n}(z;J^{(n)})}
\]
By Theorem~\ref{TA.1.2}, the right-hand side goes to zero as $n \to \infty$.

\smallskip
(iii) For $z\in\bbD$, $u_n (z;J)$ is the Weyl solution of \eqref{A.1.37},
that is, the unique up to a  multiple solution of \eqref{A.1.37} which
is $\ell^2$ at infinity (since \eqref{A.1.37a} implies $u_n\sim z^n\in
\ell^2$). Thus, by general principles, the matrix elements of $R(z)
\equiv (z+z^{-1} -J)^{-1}$ are
\[
R(z)_{nm} = \f{p_{\min (n,m)-1} (z) u_{\max(n,m)} (z)}{W(p,u)}
\]
where $W(p,u)$ is the Wronskian of $p,u$. Evaluating the $W$ at $n=0$,
\[
W(p,u) =u_0 (z;J) = u(z;J)
\]
Thus (the Weyl formula for the Weyl-Titchmarsh function)
\begin{align}
R_{11}(z) &= \f{u_1(z;J)}{u_0 (z;J)} \lb{A.1.39b} \\
&= \f{a_1^{-1} z u(z;J^{(1)})}{u(z;J)} \lb{A.1.40}
\end{align}
which is \eqref{A.1.38}.

\smallskip
(iv) Since $z^{-n} u_n (z;J)\to 1$ as $n\to\infty$, $u_0 (z;J)$ and
$u_1 (z;J)$ cannot have simultaneous zeros. Thus, by \eqref{A.1.40},
zeros of $u$ in $\bbD$ are precisely poles of $M(z)$, and the order
of the zeros is the order of the pole. Since $M$ has only simple poles
in $\bbD$ and precisely at points with $z=\beta$, $\beta + \beta^{-1}$
an eigenvalue of $J$, $u$ has only simple zeros at the prescribed points.

\smallskip
(v) If $z\in\partial\bbD$, $u_n (z;J)$ and $u_n (z^{-1};J)$ solve the same Jacobi
equation, and so their Wronskian is constant. Since $a_n\to 1$ and $u_n (z;J)\sim z^n$,
the Wronskian at infinity is $z-z^{-1}$; while at zero, it is given in terms of $u_0$ and
$u_1$, so
\begin{equation} \lb{A.1.41}
u_1 (z;J) u_0 (z^{-1};J) -u_1 (z^{-1};J) u_0 (z;J) = z-z^{-1}
\end{equation}
for $z\in\partial\bbD$. On the other hand, $u_n (z;J)$ is real if $z$ is real so, by
analyticity and continuity, $u_n (\bar z;J) =\ol{u_n (z;J)}$. Since $z^{-1}=\bar z$ for
$z\in\bbD$, \eqref{A.1.41} becomes
\begin{equation} \lb{A.1.42}
\Ima (u_1 (e^{i\theta},J)\, \ol{u_0 (e^{i\theta},J)}\,) =\sin\theta
\end{equation}

\eqref{A.1.42} implies for $\theta\neq 0$, $u_0 (e^{i\theta},J)\neq 0$ and
if $u_0 (\pm 1) =0$, then with $u'_0 (\pm 1, J)\equiv\lim_{\theta\downarrow 0}
u_0 (e^{i\theta},J)/(i\theta)$,
\[
-\Real (u_1 (\pm 1,J) u'_0(\pm 1, J)) =1
\]
proving that $u'_0\neq 0$.

\smallskip
(vi) $u_n (z^{-1};J)$ and $u_n (z;J)$ obey the same equation, so
\eqref{A.1.41} holds. This implies it cannot happen that $u_0 (z;J)$
and $u_0 (z^{-1};J)$ are both zero.

\smallskip
(vii) By \eqref{A.1.39b},
\[
\Ima M(e^{i\theta}) = \abs{u_0 (e^{i\theta},J)}^{-2}
\Ima (u_1 (e^{i\theta}, J)\, \ol{u_0 (e^{i\theta}, J)}\,)
\]
which, given \eqref{A.1.42}, implies \eqref{A.1.39a}.
\end{proof}

{\it Warning.} In \cite{Jost1}, we constructed a Jost solution $u_n (z;J)$
for certain $J$'s and $z\in\bbD$ which may not obey \eqref{A.1.21}. It
can happen that for such $J$'s that $u(z;J)$ has a zero boundary value
as $z= re^{i\theta}\to e^{i\theta}$. This is because while $z^{-n} u_n
(z;J)\to 1$ for $z\in\bbD$, this may not be true for the boundary values.

\bigskip

\end{document}